\documentclass{article}
\usepackage{graphicx}

\begin{document}

\title{\textbf {Dido’s Problem: When a myth of ancient literature became a problem of variational calculus}}
\author{Dora Musielak}
\date{}

\maketitle

\begin{abstract}
When introducing the calculus of variations, we may invoke Dido’s problem to illustrate the most fundamental variational problem: to find the curve of given perimeter which bounds the greatest area. This type of problem led mathematicians to invent solution methods of maxima and minima, and the genesis of variational calculus as a distinct branch of analysis. Dido’s problem was inspired by the mythical tale of the foundation of Carthage (ancient city in North Africa) by a Phoenician princess as told independently by Roman poet Virgil, and by Latin historian Justinus in the first two centuries B.C. Historians have debated the facts surrounding Carthage’s birth; however, contemporary mathematicians have accepted the vague events described by Virgil in his {\it Aeneid}, adding details to Dido’s story to extrapolate a few verses and use as a basis for the isoperimetric theorem. Was Leonhard Euler or Lord Kelvin who first interpreted Virgil’s poem as Dido’s problem of variational calculus? In this article I attempt to resolve a question of historical attribution to identify who first defined Dido’s problem.
\end{abstract}



\section{Introduction}

In 1937, Karl Menger{\footnote{Menger (1937)}  wrote: “The first human being to solve a problem of calculus of variations seems to have been Queen Dido of Carthage.” Contemporary mathematics books{\footnote{See for example, Brunt (2004); Freguglia and Giaquinta (2016); Coppersmith (2017); Nahin (2004); Rojo and Bloch (2018).}}  go much further than that by adding details to Dido’s story taken from Virgil’s {\it  Aeneid}, alleging that  Dido established Carthage, an ancient city in modern Tunisia, by application of the isoperimetric property of the circle to secure the largest area of land she bought upon arrival to North Africa. Here are two examples, where I use Italics to highlight details of Dido’s story that are not in Virgil’s poem:\\

{“Dido was a Carthaginian queen (ca. 850 B.C.?) who came from a dysfunctional family. Her brother, Pygmalion, murdered her husband (who was also her uncle) and Dido, with the help of various gods, fled to the shores of North Africa with Pygmalion in pursuit. Upon landing in North Africa, legend has it that she struck a deal with a local chief to procure as much land as an oxhide could contain. {\it She then selected an ox and cut its hide into very narrow strips, which she joined together to form a thread of oxhide more than two and a half miles long. Dido then used the oxhide thread and the North African sea coast to define the perimeter of her property} ... {\it it is clear that Dido sought to enclose the maximum area within her ox and the sea}. The city of Carthage was then built within the perimeter defined by the thread and the sea coast. Dido called the place Byrsa meaning hide of bull."}{\footnote{Brunt (2004), pp. 14-15.}}\\

{“Dido ... using the seashore (given as straight) as part of the boundary, {\it she laid out the hide-strip to enclose the maximum possible area, which she “knew” would be in the shape of a semicircle}."}{\footnote{Nahin (2004), p. 45.}}\\

If these accounts were based on fact, then Dido would be the first woman in humanity’s history to understand a mathematical principle, much before the first mathematicians in recorded history. Since Carthage was founded in 814 B.C., Dido was born centuries before Thales of Miletus (c. 624-548 BC), Pythagoras of Samos (c. 570-490 BC), and Euclid of Alexandria (325-265 BC), and much earlier than the Greek mathematicians who dealt with isoperimetric problems, e.g. Zenodorus (c. 200-140 BC) who wrote {\it On Isoperimetric Figures}; this work is lost but details are found in the commentaries by Theon of Alexandria (335-405 AD), and by Pappus of Alexandria (290-350 AD). In his {\it  Mathematical Collection}, Pappus presented results from ancient isoperimetry studies{\footnote{According to Pappus, the first proof of the isoperimetric property of the circle (using geometric arguments) is due to Zenodorus.}}  but he did not mention Dido.

Dido’s problem is now taught as the most fundamental isoperimetric problem: for a fixed perimeter, determine the shape of the closed, planar curve that encloses the maximum area. The answer is the circle, as any grammar school child knows, but in variational calculus the solution is determined by an analytical method introduced by Leonhard Euler and refined by Joseph-Louis Lagrange.

\section{Dido in Ancient Literature}

The story of Dido and the foundation of Carthage was immortalized by Virgil in his {\it Aeneid}, and by third century Roman historian Justinus in his {\it Epitoma historiarum Philippicarum Pompei Trogi}. However, in these stories there is absolutely no mention of Dido enclosing a circular shape for the purchased land with the string of hide, or of her using the knowledge that the circle encloses the largest area. 

Virgil wrote in the {\it  Aeneid},{\footnote{Written between 29 and 19 BC, this epic poem in 12 books tells the story of the foundation of Rome from the ashes of Troy. Virgil describes the foundation of Carthage by Dido in Book I: 297-371.}}  Book 1, lines 365-368, referring to Dido and her people arriving to Africa: “They came to this place, and bought land, where you now see the vast walls, and resurgent stronghold, of new Carthage, as much as they could enclose with the strips of hide from a single bull, and from that they called it Byrsa.”{\footnote{Line 365: {\it Devenere locos, ubi nunc ingentia cernis moenia surgentemque novae Karthaginis arcem, mercatique solum, facti de nomine Byrsam, taurino quantum possent circumdare tergo}.}} 

Justinus, who refers to Dido by her Phoenician name Elissa, wrote in Book XVIII: “By this means some respite was given to the fugitives; and Elissa, arriving in a gulf of Africa, attached the inhabitants of the coast, who rejoiced at the arrival of foreigners, and the opportunity of bartering commodities with them, to her interest. Having then bargained for a piece of ground, as much as could be covered with an ox-hide, where she might refresh her companions, wearied with their long voyage, until she could conveniently resume her progress, she directed the hide to be cut into the thinnest possible strips, and thus acquired a greater portion of ground than she had apparently demanded; whence the place had afterwards the name of Byrsa.”{\footnote{Marcus Junianus Justinus, {\it Epitoma historiarum Philippicarum Pompei Trogi} (Epitome of the Philippic History of Pompeius Trogus). Translated by Rev. John Selby Watson. London: Henry G. Bohn, Convent Garden (1853).}}

Thus, if neither Virgil nor Justin provided the details given by contemporary mathematicians about Dido’s problem as we know it, who did? how Dido’s mythical tale became Dido’s problem? Surely, only a scientist would have made the connecting leap between “bought land, …, as much as they could enclose with the strips of hide from a single bull …” (Virgil, 1st century BC), and interpreted these words as “she laid out the hide-strip to enclose the maximum possible area, which she “knew” would be in the shape of a semicircle” (Nahin, 2004).

\section{Isoperimetry and Calculus of Variations}

Isoperimetrics provided the roots for the development of the variational methodology, starting with the observation made by ancient scholars that most motion appears to be in either straight lines or circles. The definition of a straight line as the shortest path between two points was an early expression of a minimization principle known to ancient geometers. The isoperimetric problems considered in antiquity (e.g. the circle in the plane and the sphere in three-dimensional space were known as the least perimeter figures to enclose a given area and a given volume, respectively) were solved by geometric means. 

Pappus gave credit to Zenodorus (200-140 BC) for solving for the optimal form of a maximum area surface for a given perimeter. He also expounded the work of Hero (or Heron) of Alexandria (c. 10-75 AD) who studied the optics of reflection, finding that reflected light travels in a way that minimizes its travel time. The law of reflection of light—that the angle of incidence equals the angle of reflection—was well established since ancient times. In his {\it Catoptrics}, Euclid noted that light travels in straight lines and described the law of reflection (300 BC). Hero showed by a geometrical method that the actual path taken by a ray of light reflected from a plane mirror is shorter than any other reflected path that might be drawn between the source and point of observation. 

Ancient Greeks first conceived the idea that Nature selects the shortest, easiest and most direct path in moving objects between points. In the seventeenth and eighteenth centuries, ideas about the economy of Nature continued preoccupying philosophers and scientists. Finding analytic solutions to more complicated problems of maxima and minima attracted the greatest mathematicians such as Fermat, Newton, Leibniz, the Bernoulli brothers (Jacob and Johann I), Euler, Lagrange, and Maupertuis. 

Perhaps inspired by Hero’s reflected light minimization problem, Pierre de Fermat (1601-1665) showed that the time required for a light ray to traverse a neighboring virtual path differs from the time actually taken by a quantity of the second order. This is known as Fermat’s principle of least time.

Newton (1643-1727) examined the motion of bodies in a resisting medium, finding the shape of the body that renders its resistance minimal. 

In June of 1696, Johann Bernoulli (1667-1748) posed the following problem as a challenge to mathematicians: Given two points A and B in a vertical plane, find the path AMB down which a movable point particle M must, by virtue of its weight, will traverse in the shortest possible time (assumes that M’s acceleration is due only to gravity). This is the famous Brachistochrone (from the Greek {\it brachistos}, shortest, and {\it chronos}, time) problem, later also called the problem of least time descent. The brachistochrone problem does not have a trivial solution; the Bernoulli brothers (Jacob and Johann I), Newton, Leibniz and l’Hôpital solved the problem correctly, each using a different approach.{\footnote{Fregulia and Giaquinta (2016), pp. 3-4.}} 

The initial investigations in the maxima and minima principles carried out by Leonhard Euler began from a study of the work of these mathematicians, especially motivated by the work of Jacob Bernoulli and prompted by his teacher Johann Bernoulli. The latter drew his attention to a problem of geodesic lines in a letter he sent Euler in St. Petersburg in 1728, which led Euler to conceive in early 1729 an analytical method by which, on any surface, whether convex or concave, the shortest line can be drawn between two points.{\footnote{Euler (1732)}} Euler solved other isoperimetric problems, obtaining results to help him establish the analytical foundations of the calculus of variations. 

Euler invented variational calculus as a distinct branch of analysis precisely to systemize the solution methods of maxima and minima, as brilliantly introduced in his 1744 book {\it  Methodus inveniendi lineas curvas maximi minimive proprietate gaudentes, sive solutio problematis isoperimetrici lattissimo sensu accepti}, the first treatise on calculus of variations.{\footnote{Euler (1744)}}  With Euler’s approach, the calculus of variations yielded a method for finding an extremum of a quantity that is expressible as an (variational) integral. 

Euler’s {\it Methodus inveniendi} represented a substantial break with the then established tradition set for by his predecessors, including his earlier work in the subject.{\footnote{Fraser (1993)}}  In this treatise, Euler formulated the variational principle of mechanics, which is the principle of least action now attributed to Maupertuis:  For a given projected body, denote its mass by $M$, half the square of its velocity by $v$, the arclength element by $ds$. Then, among all curves passing through the same pair of endpoints, the desired curve is the one that minimizes the integral $\int M ds v^{1/2}$. Details on how Euler formulated the principle are provided by Goldstine (1980), and by Freguglia and Giaquinta (2016).{\footnote{Goldstine (1980), p. 101;  Freguglia and Giaquinta (2016), pp. 181-189.}}  

Euler remarked: “Since the structure of the universe was made most perfect as designed by the wisest Creator, nothing in the world will occur in which no maximum or minimum rule is shining forth; wherefore there is absolutely no doubt that all the effects of the world can be equally successfully determined from final causes by means of the maximum and least methods, and from the efficient causes themselves."{\footnote{Euler (1744), Additamentum I.}}

Considered as the first variational treatment of mechanics, Euler's principle of least action contributed significantly to analytic mechanics and ultimately to the fundamental underpinnings of twentieth-century physics, including general relativity and quantum mechanics.

Euler was also known as being able to recite Virgil’s {\it Aeneid} by heart. Did he interpret Dido’s tale as Dido’s isoperimetric problem? 

\section{Defining Dido’s Problem in the Calculus of Variations}

A casual survey of the history of mathematics books written in the eighteenth and nineteenth century yields no clues as to when or how Dido’s mythical story became part of variational calculus. It required a person with mathematical brilliance and fertile imagination to connect ancient myth with mathematics. Two names emerge as potential candidates: Leonhard Euler (1707-1783), the originator of the calculus of variations, and British mathematician, physicist and engineer William Thomson, known in physics as Lord Kelvin (1824-1907). 

\subsection{Leonhard Euler}

In {\it Methodus inveniendi}, Euler gives the following example to demonstrate his analytical method: to find among all admissible curves, enclosing a given area, the one of least length. Figure 1 is Euler’s sketch to demonstrate that the curved arc of a circle, BM, is minimum. In his own words:{\footnote{Euler (1744), Chapter IV, p. 135, Exemplum II: 9. {\it Super axe AP construere lineam BM, ita comparatant, ut, abscissa area ABMP datæ magnitudinis, arcus curvæ BM illi areæ respondens sit omnium minimus}.}}  \\

“On the axis AP construct the line BM, so that, when the area ABMP of a given size is cut off, the curved arc BM corresponding to that area is the minimum of all.” After solving his variational integral, Euler shows the solution curve to be an arc of a circle with center somewhere on the line AP, for example, at C in Fig. 1.\\

\begin{figure}
\begin{center}
\includegraphics{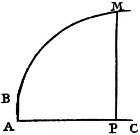}
\end{center}
\caption{Euler's sketch from {\it Methodus inveniendi} (1744)}
\label{figure}
\end{figure}

But neither Euler’s {\it  Methodus inveniendi} nor his other published memoirs in the field ever mention Dido. 

In October 1783, the Marquis de Condorcet{\footnote{Condorcet, Jean-Antoine-Nicolas de Caritat marquis de (1743-1794). }}  gave the {\it Éloge d’Euler} to the members of the Académie des Sciences in Paris. In this solemn eulogy, Condorcet expounded on Euler’s genius and suggested that a verse from the {\it Aeneid} had given Euler the first idea for a memoir on a question of Mechanics.  In Condorcet’s own words:\\

{\it L’étude de la Littérature ancienne et des Langues savantes avait fait partie de son éducation ; il en conserva le goût toute sa vie, et n’oublia rien de ce qu’il avait appris ; mais il n’eut jamais ni le tems ni le désir d’ajouter à ses premières études : il n’avait pas lu les Poètes modernes, et savait par cœur l’Eneide. Cependant M. Euler ne perdait pas de vue les Mathématiques, même lorsqu’il récitait les vers de Virgile ; tout était propre à lui rappeler cet objet presque unique de ses pensées, et on trouve dans ses ouvrages un savant Mémoire sur une question de Mécanique, dont il racontait qu’un vers de l’Eneide lui avait donné la première idée.}{\footnote{Éloge d’Euler Prononcé à l’Académie, par de Condorcet, Histoire de l’Académie royale des sciences ... 1783, p. 64.}}  [The study of ancient literature and scholarly languages had been part of his education; he retained a taste for it all his life, and forgot nothing he had learned; but he never had either the time or the desire to add to his first studies: he had not read the Modern Poets, and knew the Aeneid by heart. However, M. Euler did not lose sight of Mathematics, even when he recited the verses of Virgil; everything was likely to remind him of this almost unique object of his thoughts, and we find in his works a scholarly Memoir on a question of Mechanics, of which he said that a verse from the Aeneid had given him the first idea.]\\

Euler did take verses from the {\it Aeneid} poem to use as mottos for his competing memoirs submitted to the French Academy.{\footnote{Submissions were anonymously and the memoir identified by a motto; the author’s name enclosed in a sealed envelope was opened only for the winning memoir after the judging of the contest.}} These are summarized in Table 1. Was this to what Condorcet referred to?\\

Table 1. Euler's Memoirs and Mottos taken from Virgil's {\it Aeneid}.\\
\begin{center}
\centering 
\begin{tabular}{|c| c| c|} 
\hline 
Year & Memoir Title & Motto\\  
\hline
1753 & “On the movement of ships & {\it Tali remigio navis se} \\ 
(E. 413) & without the wind's force." & {\it tarda movebat}. \\ 
              & 7th winning memoir & Virg. Aeneid Liv. 5 \\ 
\hline
1759 & “Concernin pitching & {\it Insequitur clamorque virum} \\ 
(E. 415) & and rolling." & {\it stridorque rudentum}. \\ 
        & 9th winning memoir & Virg. Aeneid, Liv. 1 \\ 
\hline
1770 & “Moon Theory" & {\it Errantem que canit Lunam} \\ 
(E. 485) &  Prize for 1770 & Virg. Aeneid Liv. 1 \\ 
              & 10th winning memoir &  \\ 
\hline
1772 & “Improved Moon theory" & {\it Hic labor extremus, longarum} \\ 
(E. 486) & Prize for 1772 & {\it haec meta viarum hinc jam} \\ 
          & 11th winning memoir & {\it digressi, vestris appellimus oris} \\ 
          &                                  & Virg. Aeneid, Liv. 3 \\ 
\hline 
\end{tabular}
\end{center}

However, the mottos were carefully selected by Euler to match the research topic of the competition.{\footnote {For the significance of the mottos that Euler selected, see Musielak (2022).}} In addition to using Virgil’s verses, he also quoted from other ancient writers such as Marcus Tullius Cicero, Properci, and he composed other adages, asking Christian Goldbach for suggestions. Ultimately, Condorcet’s statement “{\it et on trouve dans ses ouvrages un savant Mémoire sur une question de Mécanique, dont il racontait qu’un vers de l’Eneide lui avait donné la première idée}" does not mean that Euler was inspired by Virgil to define Dido’s problem.

As a historian, I cannot rely on obituaries to extract factual data, even if written by an eminent scholar. The much younger Condorcet never met Euler, and the {\it  Éloge} he wrote, as most eulogies are, was based on hearsay, relaying on what the French academicians might have recalled about Euler’s life and work.


A contemporary biography (published in 2016) further implies that Euler solved Dido’s problem. The author refers to a copy of an eight-page manuscript (preserved in Moscow) that is said to contain Euler’s answer. Is this the manuscript that categorically would give Euler credit for connecting Dido’s story to variational calculus? Unfortunately, the manuscript in question is said to be “not in Euler’s own handwriting.” Thus, it diminishes its credibility. It is rather improbable that Euler, a prolific writer, would be the author of a manuscript inscribed by someone else.  Besides, he would have included this solution in a paper published in 1764, where Euler summarized the results of the Calculus of Variations in terms of the variational operator. 

Joseph-Louis Lagrange (1736-1813) expanded the variational calculus. In his second letter to Euler dated August 1755, Lagrange outlined his $delta$-algorithm (for solving constrained optimization problems), an approach Euler embraced, prompting him to conceive the term calculus of variations. In the abstract of a memoir published in 1764, Euler credits Lagrange for enriching the science{\footnote{Euler (1764)… {\it ex quo Auctori occasio est oblata hanc scientiam novo Calculi genere locupletandi, quem Calculum variationum appellat et cuis elementa hic tradere ac dilucide explicare constituit.} }}  Their combined work led eventually to the Euler–Lagrange equations, which are the equilibrium equations for minima of variational integrals.{\footnote {See Freguglia and Giaquinta (2016) for an excellent presentation of the Euler-Lagrange equations, including a historical perspective.}} 

Five years after Euler died, Lagrange published {\it Mécanique analytique}, his compendium on analytical mechanics, using variational ideas to present mechanics from a unified analytic viewpoint.  When teaching at the École Polytechnique in 1799, Lagrange published {\it Leçons sur le calcul des fonctions} and explained the method of variation. Lagrange provided a brief overview of the development of problems of maxima and minima, referring only to Greek mathematician Apollonius (262-190 BC), which dealt exclusively with the largest and smallest straight lines which can be drawn from given points to the arcs of conic sections.{\footnote {Lagrange (1806). {\it Les questions de maximis et minimis n’ont pas été incounues aux anciens géomètres ; car on a un livre entier d’Apollonius, qui traite presqu’uniquement des plus grandes et des plus petites lignes droites qui peuvent être menées de points donnés aux arcs des sections coniques.} p. 424.}}  Dido’s problem is not mentioned here nor in Lagrange’s other published works. 

\subsection{William Thomson, Lord Kelvin}

The first instance in which Dido’s name appear in the context of interest is found in a public lecture delivered by William Thomson in 1893. A great physicist known today as Lord Kelvin, his contributions include a major role in the development of the second law of thermodynamics; the absolute temperature scale (measured in kelvins); the dynamical theory of heat; the mathematical analysis of electricity and magnetism, including the basic ideas for the electromagnetic theory of light; and much more. He brought together disparate areas of physics—heat, thermodynamics, mechanics, hydrodynamics, magnetism, and electricity. Lord Kelvin played a key role in the final synthesis of 19th-century science, which viewed all physical change as energy-related phenomena.{\footnote{Gray (1910)}} 

Lord Kelvin related Dido’s clever approach to bargaining for land as follows, using the sketch in Fig. 2 to illustrate Dido’s problem:\\

“… She cut the ox-hide into an exceedingly long strip, and succeeded in enclosing between it and the sea a very valuable territory on which she built Carthage. In Dido’s problem the greatest value of land was to be enclosed by a line of given length. If the land is all of equal value the general solution of the problem shows that her line of ox-hide should be laid down in a circle. It shows also that if the sea is to be part of the boundary, starting, let us say, southward from any given point, A, of the coast, the inland bounding line must at its far end cut the coast line perpendicularly. Here, then, to complete our solution, we have a very curious and interesting, but not at all easy, geometrical question to answer: What must be the radius of a circular arc, ADC, of given length, and in what direction must it leave the point A, in order that it may cut a given curve, ABC, perpendicular at some unknown point, C?”{\footnote{Thomson (1894), p. 572-574.}}\\

\begin{figure}
\begin{center}
\includegraphics{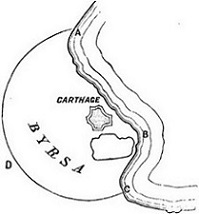}
\end{center}
\caption{Dido’s problem as described by Lord Kelvin in 1893.}
\label{figure}
\end{figure}

Lord Kelvin added that having enough mathematics knowledge, Dido would determine that the boundary had to be a circle. Of course, as illustrated in Fig. 2, she would have given the boundary a different curvature in different parts to gain as much as possible of the more valuable parts of the land offered to her, “even though difference of curvature in different parts would cause the total area enclosed to be less than it would be with a circular boundary of the same length.”{\footnote{Ibid., p. 574.}}

Today, taught as introduction to calculus of variation, the solution of Dido’s problem requires an extremization solution under constraint, that is, we maximize the area, $A = \int y dx$, subject to the condition that the arc, $L=\int ds$  is of a given length $L$. In other words, we wish to maximize the integral $A$ subject to the condition that another integral $L$ has a given constant value. Note this is an optimization problem with constraints where we use Lagrange’s strategy for finding the local maxima and minima of a function subject to equality constraints.

It is clear that, without an original reliable source, I cannot conclude that Euler defined Dido’s Problem for the first time, inspired by Virgil’s {\it Aeneid}, as Condorcet implied. The evidence points to Lord Kelvin who described the problem in 1893. And as he stated, whether severe critics will call Dido’s story mythical or allow it to be historic, it is nevertheless full of scientific interest. 

As for me, Dido’s Problem is an excellent example to introduce students to the calculus of variations, as it expresses a perfectly definite case of isoperimetrics, illustrating the fundamental principles introduced by Euler and Lagrange in the eighteenth century.

\section{Dido and Ancient Mathematics}

Nothing is known about Dido’s knowledge. Being a Phoenician princess, it is highly probable that she was well educated. What we glimpse from Virgil’s and Justinus’s tales is that Dido was a formidable woman, smart, ambitious, a foreign leader that left the city of Tyre (on the coastline of modern Lebanon) with her faithful followers, navigated the waters of the Mediterranean Sea and landed in the coast of North Africa. There, she established and ruled Carthage (modern day Tunis), an important port city that rose to the height of its power in the second century BC, before Rome became supreme and took over that region.

For the ancient cultures that flourished around the Mediterranean, geometry was fundamental to their development. The Babylonians thriving in the Mesopotamian River Valley engaged in commerce through the Mediterranean, and this required considerable mathematical skills. Clay tables preserve records of what they knew. For instance, clay tables from Babylon, located in the southern part of Mesopotamia, about fifty miles south of present-day Baghdad (Iraq) suggest that the Babylonian had an advanced knowledge of geometry and arithmetic.

In fact, some scholars believe that the Babylonians knew the Pythagorean theorem a thousand years before Pythagoras of Samos. At Susa, an ancient city over two hundred miles from Babylon, a set of tablets were discovered in 1936, which contain the ratios of areas and perimeters of regular polygons to their respective side lengths. The best known surviving tablet (estimated to be from between 1900 and 1600 BC) contains a list of Pythagorean triples. This suggests that the Babylonians had knowledge of the Pythagorean theorem, as well as certain algebraic identities.

Moreover, that women knew mathematics in ancient times has been extensively documented. For example, the Pythagorean society included women, some of which became famous such as mathematician Theano, who was married to Pythagoras. In the dedication of his {\it Introduction to Harmonics}, ancient mathematician and music theorist Nicomachus of Gerasa (c. 60-120 AD) addresses the lessons to a lady, one of his students.{\footnote{Biographical Note of Nicomachus, in Great Books of the Western World, Robert Maynard Hutchins, Ed., Vol. 11, p. 807.}}  In this book, also known as Manual of Harmonics, Nicomachus dealt with the theory of music, a version of Pythagorean harmonics, in which he assigned number and numerical ratios to notes and intervals. And of course, we know about Hypatia of Alexandria (c. 370-415 AD) considered the first woman scholar to attain eminence as mathematician and astronomer.{\footnote {Musielak (2020), pp. 206-207.}} 

I believe that Dido was educated in mathematics, and so she used the theorem of isoperimetry to outsmart the king who sold her the piece of land in the northern tip of Africa (today's Tunisia). Therefore, Dido’s Problem should be viewed not only to illustrate a fundamental problem of variational calculus but also as a lesson in the history of mathematics and the role ancient women played in its development. \\

{\bf References}\\

Brunt, B. van (2004). The Calculus of Variations. Published by Springer New York ISBN: 978-0-387-40247-5 DOI: 10.1007/b97436. 

Coppersmith, J. (2017) The Lazy Universe: An Introduction to the Principle of Least Action. Oxford University Press.

Euler, L. (1732). De linea brevissima in superficie quacunque duo quaelibet puncta iungente (On the shortest line joining two points on a surface). Commentarii academiae scientiarum Petropolitanae, Volume 3, 1732, pp. 110-124. (E. 9)

Euler, L. (1738). Problematis isoperimetrici in latissimo sensu accepti solutio generalis (On isoperimetric problems in the widest sense). Commenturii Academiae Scientiarum Petropolitanae 6 (1732/3), 123-155. Opera Omnia, 125, 13-40. (E. 27)

Euler, L. (1741). Curvarum maximi minimive proprietate gaudentium inventio nova et facilis (New and easy method of finding curves enjoying a maximal or minimal property). In Commentarii Academiae Scientiarum Petropolitanae 8 (1736). 159-190. Reprinted in Euler, L. Opera Omnia, I 25, 54-80. (E. 56)

Euler, L. (1744). Methodus inveniendi curvas h’neas maximi minimive proprietate gaudentes sive solution problematis isoperimetrici latissimo sensu accepti. Lausanne, Genf: M.-M. Bousquet. Reprinted in Euler, L. Opera Omnia, I 24. (E. 65). According to Eneström, Euler completed the manuscript of this work by April 1743.

Euler, L. (1764). Elementa calculi variationum (Elements of Calculus of Variations). Novi commentarii academiae scientiarum Petropolitanea 10 (1764), 1766, pp. 51-93. This research (E. 296) was presented at the Berlin Academy in 1756.

Freguglia, P. and Giaquinta, M. (2016). The Early Period of the Calculus of Variations. Published by Birkhäuser. 

Gelfand, I. M. and Fomin, S. V. (1963). Calculus of Variations. Revised English Edition Translated and Edited by R. A. Silverman Prentice-Hall, Inc. Englewood Cliffs, NJ.

Goldstine, H. (1980). A History of the Calculus of Variations from the 17th through the 19th Century. Springer-Verlag.

Gray, A. (1910). The Life of William Thomson, Baron Kelvin of Large. Nature 83, 61–65 (1910).

Lagrange (1806). Leçons sur le calcul des fonctions. Nouvelle édition, revue, corrigée et augmentée par l’auteur [J.-L. Lagrange]. Initially published as lecture notes in 1799 when Lagrange was teaching at the Ecole Polytechnique and reprinted in 1804. In 1806, Lagrange published a new edition containing two new lessons.

Menger, K. (1937). What is Calculus of Variations and What Are Its Applications? The Scientific Monthly 45 (3) (1937), 250-253.

Musielak, D. (2022). Leonhard Euler and the Foundations of Celestial Mechanics. Springer History of Physics Series. Springer Nature Switzerland. ISBN 978-3-031-12321-4. 

Musielak, D. (2020). Sophie Germain: Revolutionary Mathematician. Springer Biographies. Springer Nature Switzerland. ISBN 978-3030383770.

Nahin, J.P. (2004). When Least is Best. Princeton University Press, 2004.

Rojo, A. and Bloch, A. (2018). The Principle of Least Action: History and Physics. Cambridge University Press. doi:10.1017/9781139021029

Thomson, W. (1894). Popular Lectures and Addresses by Sir William Thomson (Baron Kelvin) in Three Volumes. Nature Series, MacMillan and Co. London 1894.\\

$\hskip0.7in$                              {Dora Musielak; University of Texas at Arlington, 6 January 2023}

\end{document}